%% Plain TeX
%%
%% Numerotation automatatique sous plain: Ch. Sorger
\newcount\secno
\newcount\prmno
\newif\ifnotfound
\newif\iffound

\def\namedef#1{\expandafter\def\csname #1\endcsname}
\def\nameuse#1{\csname #1\endcsname}

\long\def\ifundefined#1#2#3{\expandafter\ifx\csname
  #1\endcsname\relax#2\else#3\fi}
\def\hwrite#1#2{{\let\the=0\edef\next{\write#1{#2}}\next}}

% Working with lists
\toksdef\ta=0 \toksdef\tb=2
\long\def\leftappenditem#1\to#2{\ta={\\{#1}}\tb=\expandafter{#2}%
                                \edef#2{\the\ta\the\tb}}
\long\def\rightappenditem#1\to#2{\ta={\\{#1}}\tb=\expandafter{#2}%
                                \edef#2{\the\tb\the\ta}}

\def\lop#1\to#2{\expandafter\lopoff#1\lopoff#1#2}
\long\def\lopoff\\#1#2\lopoff#3#4{\def#4{#1}\def#3{#2}}

\def\ismember#1\of#2{\foundfalse{\let\given=#1%
    \def\\##1{\def\next{##1}%
    \ifx\next\given{\global\foundtrue}\fi}#2}}

% Les commandes
\def\section#1{\vskip1truecm
               \global\def\currenvir{section}
               \global\advance\secno by1\global\prmno=0
               {\bf \number\secno. {#1}}
               \smallskip}

\def\subsection{\global\def\currenvir{subsection}
                \global\advance\prmno by1
                \ind{ (\number\secno.\number\prmno) }}
\def\subsec{\global\def\currenvir{subsection}
                \global\advance\prmno by1
                { (\number\secno.\number\prmno)\ }}
\def\noeq{\global\advance\prmno by1
                \eqno{(\number\secno.\number\prmno)}}
\def\proclaim#1{\global\advance\prmno by 1
                {\bf #1 \the\secno.\the\prmno$.-$ }}

\long\def\th#1 \enonce#2\endth{%
   \medbreak\proclaim{#1}{\it #2}\global\def\currenvir{th}
\smallskip}
\def\bib#1{\rm #1}
\long\def\thr#1\bib#2\enonce#3\endth{%
   \medbreak{\global\advance\prmno by 1\bf
#1\the\secno.\the\prmno\ 
\bib{#2}$\!.-$ } {\it
#3}\global\def\currenvir{th}\smallskip}

\def\rem#1{\global\advance\prmno by 1
{\it #1 }\kern2pt \the\secno.\the\prmno$.-$}

% CROSS-REFERENCES

\def\isinlabellist#1\of#2{\notfoundtrue%
   {\def\given{#1}%
    \def\\##1{\def\next{##1}%
    \lop\next\to\za\lop\next\to\zb%
    \ifx\za\given{\zb\global\notfoundfalse}\fi}#2}%
    \ifnotfound{\immediate\write16%
             {Warning - [Page \the\pageno] {#1} No reference found}}%
                \fi}%
\def\ref#1{\ifx\labellist\empty{\immediate\write16
                 {Warning - No references found at all.}}
               \else{\isinlabellist{#1}\of\labellist}\fi}

\def\newlabel#1#2{\rightappenditem{\\{#1}\\{#2}}\to\labellist}
\def\labellist{}

\def\label#1{%
  \def\given{th}%
  \ifx\given\currenvir%
{\hwrite\lbl{\string\newlabel{#1}{\number\secno.\number\prmno}}}\fi%
  \def\given{section}%
  \ifx\given\currenvir%
    {\hwrite\lbl{\string\newlabel{#1}{\number\secno}}}\fi%
  \def\given{subsection}%
  \ifx\given\currenvir%
{\hwrite\lbl{\string\newlabel{#1}{\number\secno.\number\prmno}}}\fi%
  \def\given{subsubsection}%
  \ifx\given\currenvir%
  {\hwrite\lbl{\string%
\newlabel{#1}{\number\secno.\number\subsecno.\number\subsubsecno}}}\fi
\ignorespaces}

\newwrite\lbl

\def\openall{\openout\lbl=\jobname.lbl}
\def\closeall{\closeout\lbl}

\newread\testfile
\def\lookatfile#1{\openin\testfile=\jobname.#1
    \ifeof\testfile{\immediate\openout\nameuse{#1}\jobname.#1
                    \write\nameuse{#1}{}
                    \immediate\closeout\nameuse{#1}}\fi%
    \immediate\closein\testfile}%

\def\begin{\newlabel{ACM}{1.1}
\newlabel{0-reg}{1.1}
\newlabel{res}{1.3}
\newlabel{ext}{1.5}
\newlabel{gen}{2.1}
\newlabel{factor}{2.2}
\newlabel{smooth}{2.3}
\newlabel{Prym}{3.1}
\newlabel{jprym}{3.2}
\newlabel{conics}{4.1}
\newlabel{cubics}{4.2}
\newlabel{can}{4.3}
\newlabel{quartics}{4.4}
\newlabel{ab}{6.1}
\newlabel{bord}{6.2}
\newlabel{thethm}{6.3}
\newlabel{open}{6.4}
\newlabel{F+F}{6.5}
\newlabel{theta}{7.3}
\newlabel{tyurin}{8.1}
\newlabel{mos}{8.3}
\newlabel{lag}{8.4}}
\let\bye\end
\def\end{\closeall\bye}

%% Fin de la numerotation automatique

\input xypic 
\magnification 1250
\pretolerance=500 \tolerance=1000  \brokenpenalty=5000
\mathcode`A="7041 \mathcode`B="7042 \mathcode`C="7043
\mathcode`D="7044 \mathcode`E="7045 \mathcode`F="7046
\mathcode`G="7047 \mathcode`H="7048 \mathcode`I="7049
\mathcode`J="704A \mathcode`K="704B \mathcode`L="704C
\mathcode`M="704D \mathcode`N="704E \mathcode`O="704F
\mathcode`P="7050 \mathcode`Q="7051 \mathcode`R="7052
\mathcode`S="7053 \mathcode`T="7054 \mathcode`U="7055
\mathcode`V="7056 \mathcode`W="7057 \mathcode`X="7058
\mathcode`Y="7059 \mathcode`Z="705A
\def\spacedmath#1{\def\packedmath##1${\bgroup\mathsurround
=0pt##1\egroup$}
\mathsurround#1
\everymath={\packedmath}\everydisplay={\mathsurround=0pt}}
\def\nospacedmath{\mathsurround=0pt
\everymath={}\everydisplay={} } \spacedmath{2pt}
\def\qfl#1{\buildrel {#1}\over {\longrightarrow}}
\def\phfl#1#2{\normalbaselines{\baselineskip=0pt
\lineskip=10truept\lineskiplimit=1truept}\nospacedmath\smash 
{\mathop{\hbox to 8truemm{\rightarrowfill}}
\limits^{\scriptstyle#1}_{\scriptstyle#2}}}
\def\hfl#1#2{\normalbaselines{\baselineskip=0truept
\lineskip=10truept\lineskiplimit=1truept}\nospacedmath\smash
{\mathop{\hbox
to 12truemm{\rightarrowfill}}\limits^{\scriptstyle#1}_{\scriptstyle#2}}}
\def\gfl#1#2{\normalbaselines{\baselineskip=0truept
\lineskip=10truept\lineskiplimit=1truept}\nospacedmath\smash
{\mathop{\hbox to
12truemm{\leftarrowfill}}\limits^{\scriptstyle#1}_{\scriptstyle#2}}}
\def\diagramme#1{\def\normalbaselines{\baselineskip=0truept
\lineskip=4truept\lineskiplimit=1truept}   \matrix{#1}}
\def\vfl#1#2{\llap{$\scriptstyle#1$}\left\downarrow\vbox to
6truemm{}\right.\rlap{$\scriptstyle#2$}}
\def\mono{\lhook\joinrel\mathrel{\longrightarrow}}
\def\epi{\rightarrow \kern -3mm\rightarrow }
\def\iso{\vbox{\hbox to .8cm{\hfill{$\scriptstyle\sim$}\hfill}
\nointerlineskip\hbox to .8cm{{\hfill$\longrightarrow $\hfill}} }}

\def\sdir_#1^#2{\mathrel{\mathop{\kern0pt\oplus}\limits_{#1}^{#2}}}
\def\pprod_#1^#2{\raise 2pt
\hbox{$\mathrel{\scriptstyle\mathop{\kern0pt\prod}\limits_{#1}^{#2}}$}}

\font\eightrm=cmr8         \font\eighti=cmmi8
\font\eightsy=cmsy8        \font\eightbf=cmbx8
\font\eighttt=cmtt8        \font\eightit=cmti8
\font\eightsl=cmsl8        \font\sixrm=cmr6
\font\sixi=cmmi6           \font\sixsy=cmsy6
\font\sixbf=cmbx6\catcode`\@=11
\def\eightpoint{%
  \textfont0=\eightrm \scriptfont0=\sixrm \scriptscriptfont0=\fiverm
  \def\rm{\fam\z@\eightrm}%
  \textfont1=\eighti  \scriptfont1=\sixi  \scriptscriptfont1=\fivei
  \def\oldstyle{\fam\@ne\eighti}\let\old=\oldstyle
  \textfont2=\eightsy \scriptfont2=\sixsy \scriptscriptfont2=\fivesy
  \textfont\itfam=\eightit
  \def\it{\fam\itfam\eightit}%
  \textfont\slfam=\eightsl
  \def\sl{\fam\slfam\eightsl}%
  \textfont\bffam=\eightbf \scriptfont\bffam=\sixbf
  \scriptscriptfont\bffam=\fivebf
  \def\bf{\fam\bffam\eightbf}%
  \textfont\ttfam=\eighttt
  \def\tt{\fam\ttfam\eighttt}%
  \abovedisplayskip=9pt plus 3pt minus 9pt
  \belowdisplayskip=\abovedisplayskip
  \abovedisplayshortskip=0pt plus 3pt
  \belowdisplayshortskip=3pt plus 3pt 
  \smallskipamount=2pt plus 1pt minus 1pt
  \medskipamount=4pt plus 2pt minus 1pt
  \bigskipamount=9pt plus 3pt minus 3pt
  \normalbaselineskip=9pt
  \setbox\strutbox=\hbox{\vrule height7pt depth2pt width0pt}%
  \normalbaselines\rm}
\catcode`\@=12

\newcount\noteno
\noteno=0
\def\up#1{\raise 1ex\hbox{\sevenrm#1}}
\def\note#1{\global\advance\noteno by1
\footnote{\parindent0.4cm\up{\number\noteno}}
{\vtop{\eightpoint\baselineskip12pt\hsize15.5truecm\noindent
#1}}\parindent 0cm}
\font\san=cmssdc10
\def\ext{\hbox{\san \char3}}

\def\ten{^{\scriptscriptstyle \otimes 2}}
\font\scr=eusm10\font\scb=eusb10\font\sbr=eurm10
\def\en{\hbox{\scr\char69\hskip-0.5pt\sbr\char110\hskip-0.5pt\char100}}
\def\mo{\hbox{\scr\char77}}\def\mb{\overline{\hbox{\scb\char77}}}
\def\m{\overline{\mo}}

\def\pc#1{\tenrm#1\sevenrm}
\def\tx{\kern-1.5pt -}
\def\cqfd{\kern 2truemm\unskip\penalty 500\vrule height 4pt depth 0pt 
width 4pt\medbreak} 
\def\carre{\kern 2truemm\unskip\penalty 500\vrule height 4pt depth 0pt 
width 4pt} 
\def\virg{\raise
.4ex\hbox{,}}
\def\decale#1{\smallbreak\hskip 28pt\llap{#1}\kern 5pt}
\def\no{n\up{o}\kern 2pt}
\def\ind{\par\hskip 1truecm\relax}
\def\indp{\par\hskip 0.5truecm\relax}
\def\moins{\mathrel{\hbox{\vrule height 3pt depth -2pt width 6pt}}}
\def\rond{\kern 1pt{\scriptstyle\circ}\kern 1pt}

\def\Hom{\mathop{\rm Hom}\nolimits}
\def\Ext{\mathop{\rm Ext}\nolimits}

\def\im{\mathop{\rm Im}\nolimits}
\def\Ker{\mathop{\rm Ker}\nolimits}

\def\det{\mathop{\rm det}\nolimits}
\def\Pic{\mathop{\rm Pic}\nolimits}

\def\dim{\mathop{\rm dim}\nolimits}

\def\pf{\mathop{\rm pf}\nolimits}
\def\pl{{\bf P}^2}
\def\pq{{\bf P}^4}
\font\sc=eurm10
\def\zz{\hbox{\sc\char122}}

\frenchspacing
\input amssym.def
\input amssym
\catcode`\@=11
\mathchardef\dabar@"0\msafam@39
\def\lda{\mathrel{\dabar@\dabar@\dabar@\dabar@\dabar@\dabar@\dabar@
\mathchar"0\msafam@4B}}
\catcode`\@=12
\def\dfl#1{\buildrel {#1}\over {\lda}}

\vsize = 25truecm
\hsize = 16truecm
%\hoffset = -.15truecm
\voffset = -.5truecm
\parindent=0cm
\baselineskip15pt
\overfullrule=0pt

\begin
\centerline{\bf Vector bundles on the cubic threefold}
\smallskip
\smallskip \centerline{Arnaud {\pc BEAUVILLE}} 
\vskip1.2cm

{\bf Introduction}
\smallskip
\ind Let $X$ be a smooth cubic hypersurface in ${\bf P}^4$. In their
seminal paper [C-G], Clemens and Griffiths showed that the intermediate
Jacobian $J(X)$, an abelian variety defined analytically through Hodge
theory, is a fundamental  tool to understand the geometry of $X$. They
studied the Fano surface $F$ of lines contained in $X$, proving that the
Abel-Jacobi map embeds $F$ into $J(X)$ and induces an isomorphism
${\rm Alb}(F)\iso J(X)$. They were able to deduce from this the Torelli
theorem and
 the non-rationality of
$X$ (a problem which had resisted the efforts of the Italian
geometers)\note{Except the last one, these results had been obtained
independently by Tyurin [T1].}. 
\ind  Mumford
noticed that one can express $J(X)$ as a {\it Prym variety} and
thus  give an alternate proof for the non-rationality of $X$ ([C-G],
Appendix C); the other results of [C-G] 
can also be obtained via this approach [B2].
Later Clemens observed that one could use the twisted cubics
as well, giving an elegant parametrization of the theta divisor (see
(\ref{cubics}) below). 
At this point the cubic threefold could be considered as  well
understood, and the emphasis  shifted to other Fano threefolds. 
\ind On the other hand vector bundles on low-dimensional varieties
attracted much attention in the last two decades, notably because of
their relation with mathematical physics. Since 
 rank 2 vector bundles are well-known to be connected to codimension
2 cycles,  it seems obvious to try to relate them to the intermediate
Jacobian. Surprisingly this was done only recently: to my
knowledge the first attempt appears in [M-T], followed by [I-M] and [D].  
I would like to explain these results in this paper and to give some
applications.
\ind We will look at the simplest possible case, namely  stable rank 2
vector bundles  on $X$ with trivial determinant. Using the isomorphism
$\deg:H^4(X,{\bf Z})\iso {\bf Z}$ we consider the second Chern class
$c_2$ as a number, which is easily seen to be $\ge 2$. This leads us to
consider
 the moduli space $\mo$ of stable rank 2 vector bundles on $X$ with
$c_1=0\,,\ c_2=2$. It is a quasi-projective variety; according to
Maruyama it  admits a natural compactification $\m$ consisting of
classes of semi-stable sheaves with 
$c_1=0\,,\ c_2=2\,,\ c_3=0$. The main result is: 
\medskip
{\bf Theorem}$.-$
{\it The moduli space $\m$ is isomorphic to the intermediate Jacobian of
$X$ blown up along the Fano surface}.
\smallskip 
\ind (See Thm. \ref{thethm} below for a more intrinsic
formulation, and Cor. \ref{open} for the corresponding description of
$\mo$.)
\ind The theorem as stated is due to Druel [D];
the proof relies heavily on the results of [M-T] and [I-M], which itself
relies on [I].  Since these papers are rather technical I thought
it worthwile to write a simplified version, insisting on the geometric 
ideas. I have made this note  independent of  [I], [M-T] and
[I-M], but not of [D] to which I~will have to refer  for some delicate
technical points.
\ind The last sections contain some applications, in particular the
construction of a  completely integrable hamiltonian system
related to the situation.

\section{Vector bundles in $\mo$ and elliptic quintics}
\ind For the rest of the paper $X$ denotes a smooth cubic threefold
(over ${\bf C}$), and $\mo$ the moduli space of stable rank $2$ 
vector bundles on $X$ with $c_1=0$, $c_2=2$ (in this situation
``stable" simply means $H^0(X,E)=0$).
\thr Proposition
\bib [D] 
\enonce Let $[E]\in\mo$.  Then:
\indp{\rm a)} $H^1(X,E(n))=H^2(X,E(n))=0$ for
all
$n\in{\bf Z}$;\indp{\rm b)}
 $E(1)$ is spanned by its global
sections.
\endth\label{ACM}
{\it Sketch of proof} : Let $S$ be a general hyperplane
section of
$X$. By a result of  Maruyama, the restriction $E_{|S}$
satisfies $H^0(S,E(-1)_{|S})=0 $. The crucial point is the
vanishing of
$H^0(S,E_{|S})$.  This is somewhat delicate: Druel shows
that in any case the bundle $E(2)$ is spanned by its global
sections. Thus the zero set of a  general  section of $E(2)$
is a smooth  curve $C$. Using Chern classes and
cohomology computation Druel proves that if
$H^0(S,E_{|S})\not= 0$, 
$C$ is the projection in ${\bf P}^4$ (from an external
point)  of a curve of degree 14 and genus 5 in ${\bf P}^5$.
Such a curve is an extremal  Castelnuovo curve, so that we
know precisely how to describe it; from this one deduces that
$C$ cannot be contained in a smooth cubic threefold.
\ind Thus we have $H^0(S,E_{|S})=0$, so a non-zero section $s$ of
$E(1)_{|S}$ (which exists  since
$\chi(E(1)_{|S})=6$) vanishes along a finite subscheme $Z$ of $S$, of
degree $c_2(E(1)_{|S})=5$; moreover one shows that for $s$ general
$Z$ consists of 5 distinct points. 
Thus  we have an extension
$$0\rightarrow {\cal O}_S(-1) \qfl{s}
E_{|S}\longrightarrow {\cal I}_Z(1)\rightarrow 0\ ,$$where ${\cal
I}_Z$ is the ideal sheaf of $Z$ in $S$.  Since 
$H^0(S,E_{|S})=0$ $Z$ is not contained in a hyperplane; it follows easily
that $Z$ imposes independent conditions to degree $d$ hypersurfaces for
any $d\ge 2$, that is, $H^1(S,{\cal I}_Z(d))=0$. From the above exact
sequence we deduce $H^1(S,E(n)_{|S})=0$ for $n\ge 1$; since
$\chi(E_{|S})=0$ and $H^0(S,E_{|S})=H^2(S,E_{|S})=0$, we have also
$H^1(S,E_{|S})=0$, and we conclude that
$H^1(S,E(n)_{|S})=0$  for all $n$ by Serre duality. This
together with the exact sequence
$$0\rightarrow E(n-1)\longrightarrow E(n)\longrightarrow
E(n)_{|S}\rightarrow 0$$shows that the natural map
$H^i(X,E(n-1))\rightarrow  H^i(X,E(n))$ is surjective for
$i=1$ and injective for $i=2$. Since $H^1(X,E(n))$ and
$H^2(X,E(n))$ vanish for $|n|\gg 0$, assertion a) follows.
\smallskip 
\ind We have in particular $H^1(X,E)=H^2(X,E(-1))=0$; the
space $H^3(X,E(-2))$ is dual to
$H^0(X,E)$, which is $0$ by  stability.
 Thus the sheaf $F:=E(1)$ is $0$\tx {\it regular} in the
sense of Mumford, that is, $H^i(X,F(-i))=0$ for $i>0$. By an
easy induction argument using hyperplane sections, this
implies that $F$ is generated by its global sections
([M], lect. 14).\cqfd\label{0-reg}

 \medskip
\subsection Let me mention an amusing corollary; I refer
to [B3] for the details. The assertion a) above means that
$E$, viewed as an ${\cal O}_{\pq}$\tx module, is arithmetically
Cohen-Macaulay, that is, the ${\bf C}[X_0,\ldots,X_4]$\tx
module $\sdir_{n\in{\bf Z}}^{}H^0(X,E(n))$ admits a length 1
resolution by graded free modules. Using the alternate form
$\det:\ext^2E\rightarrow {\cal O}_X$ one shows that the
resolution can be chosen skew-symmetric; the outcome is: 
\th Corollary
\enonce The vector bundle $E$ admits a resolution
$$0\rightarrow {\cal O}_{\pq}(-2)^6\qfl{M}{\cal O}_{\pq}(-1)^6
\longrightarrow E\rightarrow 0\ ,$$where $M$ is a
skew-symmetric matrix with linear entries, such that $X$ is defined by 
the  equation $\pf(M)=0$.\cqfd
\endth\label{res}
\ind We will now relate our vector bundles to a certain type of curves 
on $X$ through Serre's construction.
\th Proposition
\enonce The zero locus of a non-zero section  $s\in
H^0(X,E(1))$ is a locally  complete intersection
 quintic curve $\Gamma \i\pq$ spanning $\pq$,  
 with trivial canonical
bundle and  $\dim
H^0(\Gamma ,{\cal O}_\Gamma )=1$. Conversely,
given a curve
$\Gamma \i X$ with the above properties, there exists a
vector bundle
$E$ of $\mo$ and a section $s$ of
$E(1)$ whose zero locus is
$\Gamma\ ;$ the pair $(E,s)$ is uniquely determined up to
automorphism. 
\ind   If $s$ is general enough, its zero locus is a
smooth elliptic quintic curve.
\endth  
\ind We will refer to the curves with the properties stated
in the Proposition as {\it elliptic quintics}. It is not
difficult to list all possible configurations for such
curves, but we will not need this.
\smallskip 
 {\it Proof} : Since $E$ is stable, any non-zero section
$s$ of  $E(1)$ vanishes along a
l.c.i curve $\Gamma \i X$, giving rise to an extension
$$ 0\rightarrow {\cal O}_X\qfl{s} E(1)\longrightarrow {\cal
I}_\Gamma (2)\rightarrow 0\ ,\noeq$$\label{ext}where ${\cal
I}_\Gamma
$ is the ideal sheaf of $\Gamma $ in $X$.
 From this exact sequence we get
\indp(\ref{ext}.{\it a}) $H^0(X,{\cal I}_\Gamma
(1))\cong H^0(X,E)=0$ by stability, hence $\Gamma $ spans
$\pq$; 
\indp(\ref{ext}.{\it b}) $H^1(X,{\cal I}_\Gamma )\cong
H^1(X,E(-1))=0$ by (\ref{ACM}), hence  $h^0({\cal O}_\Gamma
)=1$; 
\indp(\ref{ext}.{\it c})\kern6pt by restriction to $\Gamma$, an
isomorphism
$N_{\Gamma/X}\iso E(1)_{|\Gamma }$, which gives by the
adjunction formula
$\omega_\Gamma\cong{\cal O}_\Gamma  $.
\ind The degree of $\Gamma $ is  
$c_2(E(1))=5$. Since $E(1)$ is spanned by its global
sections (\ref{ACM}), $\Gamma $ is smooth for $s$ general.
\ind Conversely, let $\Gamma \i X$ be an
  elliptic quintic.  Since 
$\omega_\Gamma $ is trivial, we can apply Serre's
construction to $\Gamma $: there are canonical isomorphisms
$$\Ext^1({\cal I}_\Gamma ,\omega_X)\iso \Ext^2({\cal
O}_\Gamma ,\omega_X)\iso H^0(\Gamma ,\omega_\Gamma )\ ,$$so
up to isomorphism there is a unique non-trivial extension
$$ 0\rightarrow {\cal O}_X(-1)\qfl{s} E\longrightarrow {\cal
I}_\Gamma (1)\rightarrow 0\ ;$$since a generator of
 $H^0(\Gamma ,\omega_\Gamma )$ is
everywhere $\not= 0$,  the sheaf
$E$ is a rank 2 vector bundle with $c_1=0$ and $c_2=2$. We
have 
$H^0(X,E)\cong H^0(X,{\cal I}_\Gamma(1))=0$, so $E$ is
stable.\cqfd
\medskip
\rem{Remark}
The cubic $X$ contains smooth elliptic curves of degree 5
which are contained in a hyperplane, and are therefore {\it
not} elliptic quintics in our sense: in fact, any smooth
hyperplane section of $X$ contains many $4$\tx
dimensional linear systems of such curves. One can still
perform Serre's construction with these curves, but the
resulting vector bundle $E$ has a non-zero section, hence it
is no longer stable. 
\section{The Abel-Jacobi map; infinitesimal study}
\subsection Recall that $\Pic(X)=H^2(X,{\bf Z})={\bf Z}h$, where $h$ is
the class of a hyperplane section; we identify
$H^4(X,{\bf Z})$ to ${\bf Z}$ via the  isomorphism $\deg: \xi\mapsto
(\xi\cdot h)$. For each $d\in{\bf Z}$, we denote by
$J^d(X)$ the translate of $J(X)$ which parametrizes $1$\tx cycles of
degree $d$ on $X$. The variety $J^d(X)$ is non-canonically isomorphic to
$J(X)=J^0(X)$; morever the map
$\gamma\mapsto
\gamma+h^2$ provides canonical isomorphisms $J^d(X)\iso
J^{d+3}(X)$.
\ind Let $T$ be an algebraic variety and $\zz$  an element of the Chow
group  $CH^2(X\times T)$, such that the restriction
$\zz_t:=\zz_{|X\times\{t\}}$  has degree  $d$ for each $t\in T$; the
map $t\mapsto[\zz_t]$ of $T$ into $J^d(X)$  is algebraic. In particular,
for each  algebraic family $(\Gamma _t)_{t\in T}$ of degree $d$ curves
on $X$, the {\it Abel-Jacobi map} $\alpha:T\rightarrow J^d(X)$
given by $\alpha(t)= [\Gamma _t]$ is algebraic.  
\ind Let ${\cal F}$ be a coherent sheaf on $X$. Grothendieck has
constructed Chern classes $\tilde c_i({\cal F})\in
CH^i(X)$, which give back the usual Chern classes by applying the 
cycle map $CH^i(X)\rightarrow H^{2i}(X,{\bf Z})$. 
Let $d=c_2({\cal F})$; the class $\tilde c_2({\cal F})$
defines an element of $J^d(X)$, which we will denote by ${\goth
c}_2({\cal F})$. As above, any algebraic family $({\cal F}_t)_{t\in T}$
of coherent sheaves such that $c_2({\cal F}_t)=d$ for all $t\in T$ gives
rise to a morphism
${\goth c}_2:T\rightarrow J^d(X)$. We have in particular a morphism
${\goth c}_2:\mo\rightarrow J^2(X)$.\label{gen}
\smallskip 
\subsection Let ${\cal Q}$ be the Hilbert scheme of  elliptic
quintics contained in $X$. Serre's construction provides us
with a morphism $p: {\cal Q}\rightarrow \mo$, whose  fibre at a point
$[E]$ of $\mo$ is the  projective space ${\bf P}(H^0(X,E(1)))$. 
\ind The Abel-Jacobi map $\alpha:{\cal Q}\rightarrow J^2(X)$ which
associates to the curve $\Gamma $ the class of the cycle $\Gamma
-h^2$ factors as\label{factor}
$$\alpha: {\cal Q}\ \hfl{p}{}\ \mo\ \hfl{{\goth c}_2}{}\ J^2(X)\ .$$
\thr Proposition
\bib [M-T]
\enonce The moduli spaces  $\mo$ and ${\cal Q}$ are smooth; the
morphism  ${\goth c}_2:\mo\rightarrow 
J^2(X)$ is \'etale.
\endth\label{smooth}
{\it Proof} : a) Let us first
prove that ${\cal Q}$ is smooth. Let $\Gamma \in{\cal Q}$,
and let $N$ be its  normal bundle in $X$. Applying the functor
$\Hom^{}_{{\cal O}_{\pq}}(\ \ ,{\cal O}_\Gamma(-1) )$ to the
resolution (\ref{res}), and using the isomorphism $N^*\cong 
 E^*(-1)_{|\Gamma}$
(\ref{ext}.{\it c}), we get an exact sequence
$$0\rightarrow H^0(\Gamma ,N^*)\longrightarrow 
H^0(\Gamma ,{\cal O}_\Gamma )^6\qfl{M}  H^0(\Gamma,{\cal
O}_\Gamma (1))^6$$where the second arrow can be
identified with $H^0(\pq,{\cal O}_{\pq})^6\qfl{M}
H^0(\pq,{\cal O}_\pq(1))^6$. Since it is injective, we get 
$H^0(\Gamma ,N^*)=0$, hence by duality $H^1(\Gamma ,N)=0$.
This implies that ${\cal Q}$ is smooth at $[\Gamma ]$, of
dimension $h^0(N)=\deg N=\deg E(1)^{}_{|\Gamma }=10$. Since
${\cal Q}$ is a ${\bf P}^5$\tx bundle over
$\mo$, $\mo$ is smooth of dimension 5.
\smallskip \ind b) Let us prove now that $\alpha$ is smooth; here I 
follow closely [M-T]. Let $\Gamma $ be a smooth  elliptic quintic in $X$.
The tangent map 
$T_\Gamma (\alpha)$ maps $T_\Gamma ({\cal Q})=$ $H^0(\Gamma ,N)$
into
$T_{\alpha(\Gamma )}(J(X))=H^1(X,\Omega^2_X)^*$; using Serre
duality and the isomorphism $N^*\cong N(-2)$, we can view its
transpose as
$^tT(\alpha):H^1(X,\Omega^2_X)\rightarrow H^1(\Gamma
,N(-2))$. By Welters general results [We], this map  fits into a
commutative diagram:
$$\diagramme{H^0(X,{\cal O}_X(1))&\hfl{R}{}&H^1(X,\Omega^2_X)\cr
\vfl{r^{}_\Gamma }{}&&\vfl{}{^tT(\alpha)}\cr
H^0(\Gamma ,{\cal O}_{\Gamma}(1))
&\hfl{\partial}{}&H^1(\Gamma ,N(-2) ) }$$where $r^{}_\Gamma $
is the restriction map, $R$ the Griffiths residue isomorphism, and
$\partial$ the coboundary homomorphism associated to the exact
sequence of normal bundles
$$0\rightarrow N\longrightarrow N_{\Gamma
/\pq}\longrightarrow {\cal O}_\Gamma (3)\rightarrow
0$$twisted by
${\cal O}_\Gamma (-2)$.
\ind Since $r^{}_\Gamma $ and $R$ are bijective, we only have
to prove that $\partial$ is injective; it is sufficient to
prove  
$H^0(\Gamma ,N_{\Gamma/\pq}(-2))=0$  -- or, by duality, 
$H^1(\Gamma ,N_{\Gamma/\pq}^*(2))=0$. Using the
exact sequence
$$0\rightarrow N_{\Gamma/\pq}^* \longrightarrow
\Omega^1_ {\pq\,|\Gamma }\longrightarrow \omega^{}_\Gamma
\rightarrow 0$$and the vanishing of $H^1(\Gamma ,\Omega^1_
{\pq}(2)_{|\Gamma })$, this is equivalent to the surjectivity
of the map $H^0(\Gamma ,\Omega^1_
{\pq}(2)_{|\Gamma })\rightarrow H^0(\Gamma ,\omega_\Gamma
(2))$; it is enough to prove that the restriction map 
$H^0(\pq ,\Omega^1_
{\pq}(2))\rightarrow H^0(\Gamma ,\omega_\Gamma
(2))$ is onto. But this map can be identified with the 
{\it Wahl map}
$w: \ext^2H^0(\Gamma ,{\cal
O}_\Gamma (1))\rightarrow H^0(\Gamma ,\omega_\Gamma(2))$,
given by $w(L\wedge M)=L\,dM-M\,dL$ for $L,M\in H^0(\Gamma
,{\cal O}_\Gamma (1))$, which is known to be surjective in
this case ([W], Thm.~4.2).
\smallskip 
\ind c) Since $\alpha$ and $p$ are smooth, ${\goth c}_2$ is
smooth, hence \'etale because $\mo$ and $J(X)$ have the same
dimension.\cqfd 
\section{Reminder: $\bf J(X)$ as a Prym variety}
\subsection We need to recall some standard facts about Prym varieties
and  the corresponding description of $J(X)$; the reader familiar with 
this is advised to go directly to the next section. 
\ind Let $C$ be a curve of genus $g$ and $\pi :\widetilde{C}\rightarrow
C$  an \'etale (connected) double covering. We denote by $\sigma $ the
involution which exchanges the two sheets of the covering.  The abelian
variety $P=\im(1-\sigma )\i J(\widetilde{C})$, of dimension $g-1$, is
called the  {\it Prym variety} associated to $\pi $. It is equipped with
a principal polarization $\theta$ such that  the restriction to $P$ of
the principal polarization of
$J(\widetilde{C})$ equals $2\theta$.
\ind To describe the corresponding theta divisor it is convenient to
translate the situation in the Jacobian $J^{2g-2}(\widetilde{C})$
parametrizing divisor classes of degree $2g-2$. The divisor
classes $[D]$ such that $\pi _*D\equiv K_C$ fit into two
subvarieties $P_0$ and $P_1$ of $J^{2g-2}(\widetilde{C})$, both
isomorphic to $P$, characterized by
$$[D]\in P_\varepsilon\ \Longleftrightarrow\ h^0(D)\equiv \varepsilon\ \
({\rm mod.}\ 2)\, .$$
 In particular, the
subvariety $W\i \hbox{Sym}^{2g-2}(\widetilde{C})$ of effective divisors
$D$ such that $\pi _*D\in |K_C|$  has two irreducible components
$W_0$ and $W_1$, according to
the parity of $h^0(D)$. If $D=x_1+\ldots+x_{2g-2}$ belongs to the 
component $W_\varepsilon$, the divisor $\sigma
x_1+x_2+\ldots+x_{2g-2}$ belongs to $W_{1-\varepsilon}$. The map
$\hbox{Sym}^{2g-2}(\widetilde{C})\rightarrow J^{2g-2}(\widetilde{C})$
induces  morphisms $w_\varepsilon:W_\varepsilon\rightarrow
P_\varepsilon$, whose
 fibre  at a point  $[D]$ of $ P_\varepsilon$ is the
linear system $|D|$. It follows that
\indp$\bullet\ \,w_1:W_1\rightarrow P_1$ is birational;
\indp$\bullet\ \,w_0$ maps $W_0$ onto a divisor $\Theta \i P_0$,
which is a theta divisor on $P_0\cong P$; the map $W_0\rightarrow
\Theta$ is generically a ${\bf P}^1$\tx bundle.\label{Prym}

 \smallskip 
\subsection We now recall the description of the intermediate Jacobian
$J(X)$  as a Prym variety.
We pick up a general line
$\ell \i X$, a plane
$\pl\i\pq$ and project $X$ from $\ell $ onto $\pl$; we get a
morphism $f:X_\ell \rightarrow \pl$, where $X_\ell $ is
obtained by blowing up $X$ along $\ell $. The  fibre  of $f$
over a  point $p\in\pl$ is a conic, the residual intersection
of $X$ with the plane $\langle \ell ,p\rangle$; this conic is
smooth for general $p$, and becomes singular when $p$ lies in
the ``discriminant curve" $\Delta\i\pl$, which is a
plane quintic. When $\ell $ is general enough, $\Delta$ is
smooth, and every singular fibre is a rank $2$ conic, that
is, the union  of two distinct lines meeting at one point. These lines 
are parametrized by a curve $\widetilde{\Delta}$ with a double
\'etale covering
$\pi :\widetilde{\Delta}\rightarrow \Delta$; if $\pi ^{-1} (p)=\{x,y\}$,
we have $f^{-1} (p)=\ell _{x}\cup \ell _{y}$. As in (\ref{Prym}) we
denote by $\sigma $ the involution which exchanges $x$ and $y$.

\ind The Abel-Jacobi map
$\widetilde{\Delta}\rightarrow
J^1(X)$ gives rise to a morphism of abelian varieties
$\psi:J(\widetilde{\Delta})\rightarrow J(X)$. For $x\in\widetilde{\Delta}$,
the $1$\tx cycle $\ell _x+\ell _{\sigma x}$ is linearly equivalent to
$h^2-\ell $, hence is independent  of $p$; thus
$\psi$ annihilates $\pi ^*JC=\Ker(1-\sigma )$, and factors as
$$\psi:J(\widetilde{\Delta})\ \hfl{1-\sigma }{}\ P\ \hfl{\varphi}{}\ J(X)\
.$$
 If $D=\sum n_i\,x_i$ is a divisor on $\widetilde{\Delta}$ with $\pi
_*D\equiv 0$, we have $\sigma D\equiv -D$, hence
$$\varphi(2D)=\sum n_i\,\ell _{x_i}\ .$$
\ind The following result is stated in [C-G], App. C; a  proof can be
found for instance in [T2].
\label{jprym}
\thr Theorem \bib (Mumford)
\enonce $\varphi$ is an isomorphism of principally polarized abelian
varieties.\cqfd
\endth
\section{Low degree curves on the cubic}
\subsection I noticed in  [B1] that the Fano surface $F$ of lines
contained in
$X$ has a simple interpretation in terms of the Prym variety; this gives 
for instance an easy proof of the fact (already proved in [C-G]) that the
Abel-Jacobi map embeds $F$ into
$J^1(X)$ as a surface with minimal cohomology class ${\Theta ^3\over
3!}$.
 Iliev made the nice observation that the same method can be used
to study  higher degree curves as well [I]. His
construction works for general  conic bundles, but I will
specialize to the cubic threefold for simplicity. \label{conics}
\ind Conics are
not interesting: any conic $q$ determines a line $\ell_q$
(the residual intersection with $X$ of the plane spanned by
$q$) such that
$q+\ell_q=h^2$ in $CH^2(X)$; thus the image $F^2\i J^2(X)$ of the
variety of conics in $X$ by the Abel-Jacobi map is isomorphic to the
Fano surface $F\i J^1(X)$ through the isomorphism
$J^2(X)\iso J^1(X)$ given by $\xi\mapsto h^2-\xi$.

\ind The next case is twisted cubics. Let ${\cal T}$ be the variety of
twisted cubics contained in $X$; we denote by $\overline{\cal T}$ its
closure in  the Hilbert scheme of $X$. Let
$\alpha:\overline{\cal T}\rightarrow J(X)$ be the
Abel-Jacobi map $t\mapsto [t]-h^2$.
\th Proposition
\enonce The image of $\alpha:\overline{\cal T}\rightarrow
J(X)$  is a theta divisor  $\Theta \i J(X)$; its generic fibre is 
isomorphic to  ${\bf P}^2$.
\endth\label{cubics}
\ind I  learned this result from Herb Clemens some 25
years ago; it can be easily deduced from the parametrization  $\Theta
=F-F$ given in [C-G]. We will give a different proof, to illustrate
Iliev's method and also  to open the way for the case of quartic
curves.
\smallskip 
{\it Proof} : \ (\ref{cubics}.{\it a}) A twisted cubic  $t\i X$ is 
contained in a unique hyperplane section $S_t$ of $X$; we replace ${\cal
T}$ by the open and dense  subset consisting of those twisted cubic 
for which $S_t$ is smooth and does not contain $\ell $. On $S_t$ there
are $72$ nets ($=$ 2-dimensional linear systems) of twisted cubics. Thus
the Abel-Jacobi map factors as 
$\alpha: {\cal T}\rightarrow {\cal S}\rightarrow J(X)$,
where ${\cal S}$ parametrizes pairs $(S,L)$ of a smooth
 hyperplane section $S$ of $X$, not containing $\ell $, and a net of
twisted cubics $|L|$ on $S$. We may now add to ${\cal T}$ the
singular elements in each net of cubics, so that ${\cal T}$ becomes a
${\bf P}^2$\tx bundle over ${\cal S}$.

\ind Let ${\cal T}_\ell$ be the subvariety of cubics in ${\cal T}$
which intersect $\ell $ (that
is, which pass through the unique point of $\ell \cap S$); then  ${\cal
T}_\ell\rightarrow {\cal S}$ is a sub-$\!{\bf P}^1$\tx bundle of the
${\bf P}^2$\tx bundle  ${\cal T}\rightarrow {\cal S}$.
\ind Let $t\in{\cal T}_\ell $; the projection $p_\ell $ maps $t$ to
a  conic $c\in \pl$. The trace of $c$ on 
$\Delta$ is a divisor $\sum p_i\in |K_\Delta|$. The two
points $(x_i,y_i)$ of $\widetilde{\Delta}$ over a point $p_i$
correspond to two lines $(\ell_{x_i},\ell_{y_i})$ meeting $\ell $; since 
the projection
$t\rightarrow c$ is one-to-one, there is exactly one of these two lines, 
say
$\ell_{y_i}$, which intersects $t$; we put $D (t):=\sum x_i$.
Using the notation of (\ref{Prym}),  this gives a rational map
$D :{\cal T}_\ell \dasharrow W$.
\smallskip 
\indp(\ref{cubics}.{\it b}) {\it Claim} : $D $ maps
${\cal T}_\ell$ into  one component $W_\varepsilon$ of $W$; the  map
$D :{\cal T}_\ell\dasharrow W_\varepsilon$ is {\it birational} (in
particular, ${\cal T}_\ell $ is irreducible).

\ind Let  $D=\sum x_i$
be a general element of 
$W$; its push-down $\sum  \pi x_i$ on
$\Delta$ is cut down by a smooth conic $c$ transversal to $\Delta$. Let 
$S_c$ be the pull-back to $c$ of the conic bundle $f$:
$$\diagramme{S_c &\kern-5pt\mono \kern-5pt& X_\ell &\cr
\vfl{}{} && \vfl{}{\!f }& \cr
c &\kern-5pt\mono\kern-5pt & \pl&.}$$It is a conic bundle over $c$ with
$10$ singular fibres $f^{-1} (\pi x_i)=\ell _{x_i}\cup \ell _{\sigma
x_i}$. Let $b:S_c\rightarrow S_{{D}}$ be the blowing
down of
$\ell _{x^{}_1},\ldots,\ell _{x^{}_{10}}$; $S_{{D}}$ is a 
${\bf P}^1$\tx bundle over
$c$, thus isomorphic to one of the ruled surfaces ${\bf F}_n$ $(n\ge 0)$.
Since these surfaces have different topological type according to the
parity of $n$, we see that this parity is constant when ${D}$ varies in
each component
$W_\varepsilon$. On the other hand, let us observe that {\it the parity
of
$n$ changes with}
$\varepsilon$: consider the divisor ${D}'=\sigma x_1+
x_2+\ldots+ x_{10}$, which belongs to $W_{1-\varepsilon}$
(\ref{Prym}); the surface $S_{{D}'}$ is obtained from $S_{{D}}\cong
{\bf F}_n$ by performing an elementary transformation on the fibre at
$\pi (x_1)$, thus is isomorphic to ${\bf F}_{n\pm 1}$.

 \ind Let us now specialize to the case  where ${D}$ is the divisor 
$D(t)$ associated to a general twisted cubic $t\in {\cal T}_\ell $. Then
$t$ lies in
$S_c$, and projects down isomorphically to a curve $\bar t$ in
$S_{{D}}$. Let $E$ be the exceptional divisor of $X_\ell $. Since the 
canonical bundle of $X_\ell $ is $f^*{\cal O}_{\pl}(-2)\,(-E)\cong {\cal
O}_{X_\ell }(-S_c-E)$,  the adjunction formula gives
$K_{S_c}\equiv -E_{|S_c}$.  Therefore $-K_{S_c}\cdot t=E\cdot
t=1$, hence  by the adjunction
formula  $t^2=-1$ in $S_c$, and therefore also $\bar t^2=-1$ in
$S_{{D}}$.  But for $n>0$ the surface ${\bf F}_n$ contains a unique curve 
of negative square, namely the {\it exceptional section} of square $-n$. 
We conclude that $S_{{D}}$ is isomorphic to ${\bf F}_1$, and that
$\bar t$ is  its exceptional section; moreover the points  of $S_{{D}}$ 
which are blown up by $b$ do not lie on $\bar t$.
\ind Thus all components of ${\cal T}_\ell $ are mapped into the same
component $W_\varepsilon$ of $W$. For a general
${D}$ in $W_\varepsilon$, the surface $S_{{D}}$ is isomorphic to ${\bf
F}_1$ (a small  deformation of ${\bf F}_1$ is again isomorphic to ${\bf
F}_1$), and the exceptional section $\bar t$ does not pass through the
points which are blown up. Let
$t$ be its proper transform in $S_c$. We have $t^2=-1$, hence  $t\cdot
E=1$; since $t$ projects isomorphically onto
$c$, it is a twisted cubic  intersecting $\ell $ at one point. This 
provides (birationally) the inverse mapping of $D$. \par
\smallskip 
\indp(\ref{cubics}.{\it c})  {\it Claim} : We have a commutative diagram
$$\diagramme{
{\cal T}_\ell  & \kern-10pt\dfl{D}{}\kern-20pt  &W_\varepsilon\cr
    \vfl{\alpha}{}  &  &   \vfl{}{w_\varepsilon}  \cr
J(X) & \kern-6pt\gfl{\sim}{\varphi}\kern-6pt &P\cong
P_\varepsilon\quad\ \cr }$$where the isomorphism $P\cong
P_\varepsilon$ is the translation by an appropriate element 
$\xi\in P_\varepsilon$. 
\ind  Let $t\in{\cal T}_\ell $, and $D(t)=\sum x_i$ the corresponding
divisor on $\widetilde{\Delta}$. The commutativity of the diagram means
that the element $u(t):=\alpha(t) -\varphi
\bigl(\sum x_i -\xi\bigr)$ of $J(X)$ is zero.
 Since $\xi$ is arbitrary and ${\cal T}_\ell $ is irreducible, it
suffices to prove that
$2u(t)$ is constant when $t$ varies. In
view of the definition of
$\varphi $ (\ref{jprym}), this means that {\it the class of
$2t-\sum\ell _{x_i}$ in $J^{-4}(X)$ is independent of} $t$.

\ind Let us denote by
$f$ the class of a fibre of the ruling in
$\Pic(S_{{D}})$ and also in $\Pic(S_c)$. We have $-K_{S_{{D}}}\equiv
2\bar t+kf$ for some integer $k\in{\bf Z}$, hence 
$$ 2t-\sum \ell _{x_i}+kf\equiv -K_{S_c}\equiv E_{|S_c}\ .$$ In
$CH^2(X_\ell )$ the classes $f=f^*[p]$ and $E\cdot S_c=E\cdot
f^*[c]$ are independent of $t$, hence so is $2t-\sum \ell _{x_i}$;
pushing down in $CH^2(X)$ gives our assertion.
\smallskip 
\indp(\ref{cubics}.{\it d}) The commutativity of the diagram implies that
the generic fibre of $\alpha$ is birationally isomorphic to that of
$w_\varepsilon$. We know that $\alpha$ has positive-dimensional fibres
(\ref{cubics}.{\it a}), hence so does $w_\varepsilon$; therefore
$\varepsilon=0$. Now  the image of $w_0$ is the canonical
theta divisor in $P_0$, and its general fibre is isomorphic to ${\bf
P}^1$ (\ref{Prym}); thus the image of $\alpha$ is a theta
divisor $\Theta\i J(X)$, and its general fibre is isomorphic to  ${\bf
P}^1$. In  the factorization $\alpha: {\cal T}_\ell \rightarrow {\cal
S}\qfl{\rho} J(X)$, we see that $\rho$ is birational; therefore ${\cal
T}$ is birationally a ${\bf P}^2$\tx bundle over $\Theta$.\cqfd
  
\rem{Remark} The divisor $\Theta $  is the {\it canonical theta divisor}
of $J(X)$: its unique singular point is $0$ (see [B2]). The corresponding
divisor $\varphi^*\Theta \i P$ is the translate of the canonical theta
divisor $w_0(W_0)\i P_0$ by the element  $\pi ^*{\cal O}_\Delta(1)$
of $P_0$.\label{can}
\medskip
\ind  We now pass to  quartic curves. Let
${\cal Q}_\ell
$ be the variety of  rational normal quartic curves which are contained in
$X$ and intersect $\ell $ transversally in $2$ points; let 
$\overline{\cal Q}_\ell
$ be its closure in  the Hilbert scheme of $X$.
\th Proposition
\enonce The Abel-Jacobi map $\overline{\cal Q}_\ell \rightarrow
J^4(X)$ is surjective; its general fibre is isomorphic to ${\bf P}^1$.
\endth\label{quartics}
{\it Proof} : The proof follows closely the proof of (\ref{cubics}). 
Let $q$ be a general element of $ {\cal Q}_\ell $. Again $q$ projects
from $\ell $ onto a  smooth conic
$c\i\pl$, and defines a section of the ruled surface $S_c$. Let $\sum
p_i$ be the divisor cut down by  $c$ on $\Delta$;  we write as above
$f^{-1} (p_i)=\ell _{x_i}\cup \ell _{y_i}$, with $\ell _{x_i}\cap
q=\varnothing$, and $D(q)=\sum x_i$.
\ind As in (\ref{cubics}.{\it b}) we consider the image $\bar q$
of $q$ by $b:S_c\rightarrow S_{D(q)}$; since $E\cdot q=2$ we
find $\bar q^2=q^2=0$. This implies:
\indp1) $S_{D(q)}$ is isomorphic to
${\bf F}_0$, so by (\ref{cubics}.{\it b}) $D(q)$ lives in $W_1$; 
\indp2) $\bar q$ is one of the lines  of the
horizontal ruling.
\ind  Conversely, any general line of this ruling lifts to a quartic 
$q$ of ${\cal Q}_\ell $: thus  $D$ maps
${\cal Q}_\ell$ into  $W_1$, it is dominant, and its general fibre is
a rational curve (parametrizing the horizontal ruling).

\ind As in (\ref{cubics}.{\it c}) we get a commutative diagram
$$\diagramme{
{\cal Q}_\ell  & \kern-10pt\dfl{D}{}\kern-20pt  &W_1&\cr
    \vfl{\alpha}{}  &  &   \vfl{}{w_1}  &\cr
J(X) & \kern-6pt\gfl{\sim}{\varphi}\kern-6pt &P\cong
P_1\quad\ &\kern-25pt .\cr }$$\ind Now $w_1$ is birational, so $\alpha$
 is birationally equivalent to $D$; the Proposition follows.\cqfd 

\font\grath=eufb10
\def\cg{\hbox{\grath\char99}}

\section{$\cg_{\bf 2}$ is of degree one}
\ind We now turn back to our moduli space $\mo$ and the map ${\goth
c}_2:\mo\rightarrow J^2(X)$. We keep our general line $\ell \i X$.
\th Lemma
\enonce Let $[E]\in\mo$. The subspace of sections in $H^0(X,E(1))$
which vanish along $\ell $ has dimension $2$.
\endth 
{\it Proof} : We use the exact sequence (\ref{res}):
$$0\rightarrow {\cal O}_{\pq}(-1)^6\longrightarrow {\cal O}_{\pq}^6
\longrightarrow E(1)\rightarrow 0\ .$$
Restricting to $\ell $ we get an exact sequence
${\cal O}_\ell (-1)^6\longrightarrow {\cal O}_\ell ^6
\qfl{p} E(1)_{|\ell }\rightarrow 0\ $. The kernel $K$ of $p$ is a
rank $4$ vector bundle on $\ell $, of determinant ${\cal O}_\ell
(-2)$, with a surjective map ${\cal O}_\ell (-1)^6
\epi K$; this implies  $H^1(\ell ,K)=0$, hence $\dim H^0(\ell ,K)=2$ by
Riemann-Roch; but $H^0(\ell ,K)$ can be identified with the kernel
of the restriction map $H^0(X,E(1))\rightarrow H^0(\ell ,E(1)_{|\ell
})$.\cqfd 
\th Proposition 
\enonce  ${\goth c}_2$ induces an isomorphism of $\mo$ onto an open
subset of $J^2(X)$.
\endth
\ind We will describe precisely this open subset later (Cor.
\ref{open}).\par {\it Proof} : Recall that the variety ${\cal Q}$ of
elliptic quintics in $X$ admits a fibration $p:{\cal Q}\rightarrow \mo$,
such that the  fibre of $p$ at $[E]\in\mo$ is identified with ${\bf
P}(H^0(X,E(1)))$. Let ${\cal Q}'_\ell $ be the subvariety of ${\cal Q}$
parametrizing  elliptic quintics containing $\ell $; by the lemma the
restriction of $p$ to ${\cal Q}'_\ell $ is a sub-$\!{\bf P}^1$\tx bundle
of ${\cal Q}\rightarrow
\mo$. A curve in ${\cal Q}'_\ell $ is the union of $\ell $ and a
quartic curve, so we can view ${\cal Q}'_\ell$ as a variety parametrizing
quartic curves in
$X$, containing ${\cal Q}_\ell $; the Abel-Jacobi map ${\cal
Q}\rightarrow J^2(X)$ induces on ${\cal Q}'_\ell$, up to a constant, the
Abel-Jacobi map for this family of quartics. Since
${\cal Q}'_\ell
$ is smooth and 6-dimensional (\ref{smooth}), an easy count of constants
shows that ${\cal Q}'_\ell $ is contained in $\overline{\cal Q}_\ell
$. In the diagram
$$\xymatrix@H=14pt { 
{\cal Q}'_\ell   \ar@{_{(}->}[d]\ar[dr]_p\ar[drr]^{\alpha}\\
{\cal Q}\ar[r] & \mo \ar[r]^<<<<{{\goth c}_2}&J^2(X)}$$
we know that ${\goth c}_2$ is \'etale (Prop. \ref{smooth}) and that  
the fibres of $\alpha$ are connected (Prop.
\ref{quartics}), so ${\goth c}_2$ has degree one,
hence is an open embedding.\cqfd
\section{The boundary of $\mb$}
\subsection To complete the description of the moduli space $\m$, we
must describe the  sheaves in $\m$ which are not locally free.
Let us first describe examples of such sheaves (the systematic reference
for this section is [D]):
\ind a) Let $c$ be a smooth conic in $X$, $L$ the positive
generator of $\Pic(c)$ (so that ${\cal O}_X(1)_{|c}\cong
L\ten$). Let $E$ be the kernel of the canonical map
$H^0(c,L)\otimes_{\bf C}{\cal O}_X\rightarrow L$. Then $E$ is
a torsion free sheaf, with $c_1(E)=c_3(E)=0$ and ${\goth
c}_2(E)=[c]$. It is not difficult to check that $E$ is {\it
stable}. 
\ind b) Let $\ell,\ell '$ be two lines in $X$ (possibly
equal), and let ${\cal I}_\ell,{\cal I}_{\ell'} $ be their
 ideal sheaves. The sheaf ${\cal I}_\ell\oplus{\cal I}_{\ell'}
$ is a torsion free sheaf with
$c_1(E)=c_3(E)=0$ and ${\goth c}_2(E)=[\ell ]+[\ell ']$. It
 is clearly semi-stable but not stable. 
\smallskip 
\ind Sheaves of type a) are parametrized by  smooth conics
of $X$, hence they form a codimension $1$ subvariety $A$ of $\m$.
Sheaves of type b) are parametrized by the symmetric
square of $F$, hence they form an irreducible divisor $B\i\m$.\label{ab}
\thr Proposition 
\bib [D]
\enonce {\rm a)} $\m\moins\mo=A\cup B$.
\ind {\rm b)}  The moduli space $\m$ is smooth
and connected.\endth\label{bord}
\ind The proof of a) requires a thorough analysis of the
relationship between a sheaf in  $\m$ and its bidual; we can only refer
the reader to [D]. 
\ind Let us sketch the proof of b).
 We already know that $\mo$ is smooth
and connected; it remains to prove  that $\m$ is smooth of dimension $5$
at each point of the boundary. This is easy for a point  of $A$,
which corresponds to a stable sheaf $E$ of type a): the exact sequence
$$0\rightarrow E\rightarrow {\cal O}_X^2\rightarrow L\rightarrow
0$$gives an inclusion of $\Ext^2_X(E,E)$ into $\Ext^3_X(L,E)$; by Serre
duality, this last space is dual to $\Hom_X(E,L(-2))$, which is easily seen
to be zero. Thus $\Ext^2_X(E,E)=0$ and $\dim\Ext^1_X(E,E)=5$ by
Riemann-Roch, hence our assertion.
\ind Things become more complicated for a sheaf $E={\cal
I}_\ell\oplus{\cal I}_{\ell'}$, which is not stable. Using the
presentation of
$\m$  as a GIT quotient of a Quot scheme by a reductive
group, one shows that it is locally analytically  isomorphic to the
quotient of $\Ext^1_X(E,E)$ by the automorphism group of $E$. This group
is isomorphic to ${\bf C}^*\times {\bf C}^*$ if $\ell \not= \ell '$ and to 
$PGL(2,{\bf C})$ if $\ell =\ell '$; it is not difficult to describe its 
action on  
$\Ext^1_X(E,E)$ and to check that the quotient is smooth.\cqfd
\ind Recall that we denote by $F^2\i
J^2(X)$ the image of the variety of conics; it is isomorphic to the Fano
surface $F\i J^1(X)$ (\ref{conics}). The morphism ${\goth c}_2$  maps
$B$ onto the divisor $F+F$ in $J^2(X)$ and thus contracts only one
irreducible divisor, namely  $A$ which is mapped onto the smooth surface
$F^2$. By [Mo], Thm. 1, this implies:
\thr Theorem
\bib [D]
\enonce The map ${\goth c}_2: \m\rightarrow J^2(X)$ is isomorphic to
the blowing up of $J^2(X)$ along the Fano surface $F^2$.\carre 
\endth\label{thethm}
\ind Any conic in $X$ is linearly equivalent to a sum of
two (incident) lines, so $F^2$ is contained in $F+F$;  from
(\ref{bord} a) and (\ref{thethm}) we conclude:
\th Corollary
\enonce ${\goth c}_2$ induces an isomorphism of $\mo$ onto the
complement of the divisor $F+F$ in $J^2(X)$.\cqfd
\endth\label{open}
\rem{Remark} Let
$s:\hbox{Sym}^2F\rightarrow F+F$ be the sum map, and let
$\theta\in H^2(J(X),{\bf Z})$ be the principal polarization. Since the
class of $F$ in $H^6(J(X),{\bf Z})$ is
${\theta^3\over 3!}$,  a standard Pontrjagin product computation 
gives  $(\deg s).[F+F]=3\theta$. But $s$ can be identified with the
restriction of ${\goth c}_2$ to $B$, hence it is generically one-to-one by
the theorem; thus we find that {\it the divisor $F+F$ is algebraically
equivalent to $3\theta$}.\label{F+F}
\ind In particular,  $F+F$  is an ample divisor; thus:
\th Corollary
\enonce The variety $\mo$ is affine.\cqfd
\endth
\section{Applications}
 \ind Our first application is a Torelli-type theorem:
\th Proposition
\enonce Let $X'$ be another cubic threefold, and $\mo'\i\m'$ the
corresponding moduli spaces. If $\mo$ and $\mo'$, or $\m$ and $\m'$,
are isomorphic, then $X$ and $X'$ are isomorphic.
\endth
{\it Proof} : Any isomorphism $\mo\iso\mo'$ extends to an isomorphism
$J^2(X)\iso J^2(X')$ which maps the divisor $F+F$ onto $F'+F'$. By
Remark \ref{F+F}, the principally polarized abelian varieties $J(X)$ and
$J(X')$ are  isomorphic, hence by the Torelli theorem $X$ and $X'$ are
isomorphic. 
\ind Similarly, an isomorphism $\m\iso\m'$ induces an isomorphism
$J^2(X)\iso J^2(X')$  mapping $F^2$ onto $F'^2$, from which we
deduce  an isomorphism
$J(X)\iso J(X')$  mapping $F-F=\Theta $ onto $F'-F'=\Theta '$; again  
we conclude that $X$ and $X'$ are isomorphic.\cqfd

\subsection Once we have a parametrization of the intermediate Jacobian,
it is natural to ask if it can be used to express the theta divisor. 
Giving a theta divisor  in $J^2(X)$ is equivalent to choosing a point of
$J^2(X)$ (namely the singular point of the divisor); it is therefore 
natural to fix a smooth conic $c\i X$. Let $\Theta _c\i J^2(X)$ be the
translate by $c$ of the canonical theta divisor $\Theta \i J(X)$
(\ref{can}).
\th Proposition
\enonce Let  $\mo_c\i\mo$ be the  locus of vector bundles $E$
such that 
$E(1)$ has a non-zero section vanishing along $c$. Then the closure
$\overline{\mo}_c$ of
$\mo_c$ in $\m$ is the proper transform of  $\Theta _c$.\label{theta}
\endth
{\it Proof} : We first 
observe that  $h^0(c,E(1)_{|c})=6$ for all $[E]\in\mo$: by an easy count
of constants, a general section $s$ of
$E(1)$ does not vanish on $c$; then the exact sequence (\ref{ext})
gives by restriction to $c$ an exact sequence $0\rightarrow {\cal
O}_c\rightarrow E(1)_{|c} \rightarrow {\cal O}_c(2)\rightarrow 0$, hence
our assertion.
\ind Thus the elements $[E]$ of $\mo_c$ are those for which the
restriction map $H^0(X,E(1))\rightarrow h^0(c,E(1)_{|c})$ is not
bijective. By a well-known construction this condition defines a  {\it
divisor} in $\mo$.
\ind We now consider the factorization 
$\alpha: {\cal Q}\qfl{p} \mo\mono J^2(X)$ (\ref{factor}). Let
${\cal Q}_c$ be the subvariety of ${\cal Q}$ consisting of elliptic 
quintics containing
$c$; we have by definition $p({\cal Q}_c)=\mo_c$. A curve in ${\cal
Q}_c$ is the union of $c$ and a  curve $t$ of degree 3. Let ${\cal T}_c$
be the variety of twisted cubics in $X$ meeting $c$ transversally at 2
points; we identify it to a subvariety of
${\cal Q}_c$ by mapping $t$ to $c\cup t$.  It is 4-dimensional,  and  a
simple count of constants shows that ${\cal Q}_c\moins {\cal T}_c$ has
dimension $\le 3$. It follows that
$\mo_c=\overline{p({\cal T}_c)}$ in $\mo$. 
\ind Now by (\ref{cubics}) $\alpha({\cal
T}_c)$ is contained in the  divisor $\Theta _c \i J^2(X)$,
and therefore the divisor $\mo_c$ is contained in the restriction of
$\Theta _c$ to $\mo$; since the latter is  irreducible, they coincide, 
and
$\overline{\mo}_c$ maps onto $\Theta _c$ in $J^2(X)$. 
\ind To conclude we just have to check that $\overline{\mo}_c$ does not
contain the exceptional divisor, that is, that a general  sheaf $E$ of
type a) (\ref{ab}) satisfies $H^0(X,{\cal I}_cE(1))=0$, where  ${\cal
I}_c$ is the ideal sheaf of $c$ in $X$. Let us pick up a conic $d\i X$ 
such that the plane spanned by $d$ does not meet $c$; let $L$ be the
positive generator of $\Pic(d)$,  and let $E$ be the kernel of the
natural map
$H^0(d,L)\otimes_{\bf C}{\cal O}_X\rightarrow L$. By tensor product
with
${\cal I}_c(1)$ we get an exact sequence
$$0\rightarrow {\cal I}_cE(1)\rightarrow H^0(d,L)\otimes_{\bf C}{\cal
I}_c(1)\rightarrow L(1)\rightarrow 0\ .$$Let $V$ be the image
of $H^0(X,{\cal I}_c(1))$ in  $H^0(d,{\cal O}_d(1))$; by hypothesis it is
a base-point free 2-dimensional subspace, and
$H^0(X,{\cal I}_cE(1))$ can be identified with the kernel of the 
canonical map 
$H^0(d,L)\otimes V\longrightarrow H^0(d, L(1))  $;  the Koszul
exact sequence
$$0\rightarrow  L(-1)\rightarrow V\otimes_{\bf C}L\rightarrow
L(1)\rightarrow 0$$
gives  isomorphisms $H^0(X,{\cal I}_cE(1))\cong H^0(d,L(-1))=0$.\cqfd

\section{Application: a completely integrable system}
\subsection In the introduction of [T3], Tyurin made the following
beautiful observation. Let $X$ be a Fano threefold, and $S$ a smooth
anticanonical divisor in $X$ (so that $S$ is a K3 surface). Let
$\mo_X$ be the moduli space of stable  vector bundles on $X$ with
fixed rank $r$ and
 Chern classes  $c_1,c_2,c_3$,
and let $\mo_S$ be the moduli space of stable  vector bundles on $S$
with  rank $r$ and
 Chern classes  ${c_1}^{}_{|S}$,${c_2}^{}_{|S}$. According to Mukai
[Mu], ${\cal M}_S$ is smooth and carries a symplectic structure: at a
point $[F]\in\mo_S$, the symplectic form on
$T_F(\mo_S)=H^1(S,\en(F))$ is given by the cup-product 
$$H^1(S,\en(F))\otimes H^1(S,\en(F))\longrightarrow
H^2(S,\en(F)) \qfl{\rm Tr} H^2(S,{\cal O}_S)\iso {\bf C}\ .$$ 
\ind Let $[E]\in\mo_X$; assume that
\indp{\it a}) The restriction of $E$ to $S$ is stable;
\indp{\it b}) $H^2(X,\en(E))=0$.\par
By {\it a)} the restriction map $r:\mo_X\rightarrow \mo_S$ is defined in
a neighborhood of [E]. Using {\it b}) and Serre duality,
the exact sequence $0\rightarrow
\en(E)\otimes K_X\rightarrow \en(E)\rightarrow \en(E)^{}_{|S}\rightarrow 
0$ gives rise to  a cohomology  exact sequence 
$$ 0\rightarrow H^1(X,\en(E))\qfl{r_*} H^1(S,\en(E^{}_{|S}))\qfl{r^*} 
H^1(X,\en(E))^*\rightarrow 0\ ,$$ where $r_*$ is the tangent map to $r$
at $[E]$ and $r^*$ its transpose with respect to 
the symplectic form. Therefore  $r_*$ is injective
and  its image is a maximal isotropic subspace; thus 
$r$ {\it induces an isomorphism of  some open neighborhood of $[E]$ in
$\mo_X$ onto a Lagrangian submanifold of}
$\mo_S$.\label{tyurin}\smallskip 
\subsection  In practice, condition {\it b}) is usually difficult to 
check.  We now turn back to our situation, taking for $\mo_X$ our moduli
space
$\mo$ of rank 2 vector bundles on the cubic threefold $X$ with
$c_1=0$, $c_2=2$. We will show using the previous results
that {\it $\mo_X$ embeds as a Lagrangian submanifold of the moduli
space}
$\mo_S$. We will even show that this submanifold varies in a {\it
Lagrangian fibration} (equivalently, a completely integrable hamiltonian
system) defined on an open subset of $\mo_S$.
\ind  We pick  a   quadric $Q\i\pq$ such that
$S=X\cap Q$ is a smooth K3 surface with $\Pic(S)={\bf Z}$ [De]. We
denote by
$\mo_S$ be  the moduli space of stable rank 2 vector bundles on $S$
with  $c_1=0$,
$c_2=4$.

\def\mos{\mo_S^{\rm o}}
\th Lemma
\enonce The vector bundles $F$ in $\mo_S$ admitting a resolution
$$0\rightarrow {\cal O}_{Q}(-2)^6\qfl{M}{\cal O}_{Q}(-1)^6
\qfl{p} F\rightarrow 0\ ,$$where $M$ is a
skew-symmetric matrix with linear entries, form an open subset $\mos$ of
$\mo_S$. For such a vector bundle, the matrix $M$ is uniquely
determined up to a transformation $M\mapsto AM\,^t\!A$, with $A\in
GL(6,{\bf C})$.
\endth\label{mos}
{\it Proof} : If $F$ admits the above resolution, we have:
\indp a) For $j=0$ or $1$, the map $H^0(S,{\cal O}_S(j))^6\rightarrow
H^0(S,F(j+1))$ deduced from $p$ is surjective.\indp b) $H^1(S,F(j))=0$
for $j=0$ or $1$.
\ind Let $N$ be the kernel of the natural map
$H^0(S,F(1))\otimes_{\bf C}{\cal O}_Q\rightarrow F(1)$. In our case it is
isomorphic to
${\cal O}_Q(-1)^6$, so: 
\indp c) The sheaf $N(1)$ is spanned by its global sections. 
\ind Conditions a), b) and c) define an open subset $\mos$ of $\mo_S$.
Let $[F]\in \mos$; by a), b) and Riemann-Roch we find
$h^0(S,F(1))=h^0(Q,N(1))=6$;  hence $N$ is a rank 6
vector bundle on $Q$, and the surjective map
${\cal O}_Q^6\rightarrow N(1)$  must be an isomorphism. Thus
 we have an exact sequence $0\rightarrow {\cal O}_{Q}(-2)^6\qfl{M}
{\cal O}_{Q}(-1)^6\longrightarrow  F\rightarrow 0$. Reasoning as in 
[B3], Theorem B, we can choose bases of ${\cal O}_{Q}(-2)^6$ and 
${\cal O}_{Q}(-1)^6$ such that $M$ is skew-symmetric; then the matrix
$M$ is uniquely determined up to the action of $GL(6,{\bf C})$ ([B3],
2.3).\cqfd

\ind  We now fix our K3 surface $S$, but allow $X$ to vary in the 
 linear system $\Pi\ (\cong {\bf P}^5)$ of cubics containing $S$. For
$[F]\in\mos$, the cubic  $\pf M=0$ is well determined by $[F]$,
and belongs to
$\Pi$; we  get in this way a morphism $H:\mos\rightarrow \Pi$.   
 
\th Proposition  
\enonce $H$ is a Lagrangian fibration\note{This means, by definition, 
that the smooth part of each fibre of $H$ is a Lagrangian
submanifold.}; the fibre of $H$ at a smooth cubic
$X\in \Pi$ is isomorphic to the moduli space $\mo_X$ {\rm (}through 
the restriction map $r:\mo_X\rightarrow \mo_S\,).$
 In particular,  $r$  induces an isomorphism
of $\mo_X$ onto a Lagrangian submanifold of $\mo_S$.
\endth\label{lag}
{\it Proof} : Let  $[E]\in \mo_X$. From the exact sequence
$0\rightarrow E(-1)\rightarrow E\rightarrow E_{|S}\rightarrow 0$  and
the vanishing of $H^1(X,E(-1))$ (Prop. \ref{ACM}), we get
$H^0(S,E_{|S})=0$; since $\Pic(S)={\bf Z}$ this means that 
$E_{|S}$ is stable. 
\ind Next let us compute $H^2(X,\en(E))$. Choosing a general
section of $E$, we deduce from  the exact sequence (\ref{ext}), after
tensor product with $E\cong E^*$, an exact sequence
$$ 0\rightarrow E(-1)\longrightarrow \en(E)\longrightarrow
E(1)\rightarrow E(1)_{|\Gamma}  \rightarrow 0\ .$$Because $H^i(X,
E(j))=0$ for
$i=1,2$, this gives an isomorphism 
$H^2(X,\en(E))\iso$  $H^1(\Gamma ,E(1)_{|\Gamma })$.
Since
$E(1)_{|\Gamma }$ is isomorphic to  
$N_{\Gamma/X}$ (\ref{ext}.{\it c}) and $H^1(\Gamma
,N_{\Gamma/X})=0$ (Prop. \ref{smooth} a)), our assertion follows.
\ind  Put $F=E_{|S}$. The resolution 
$0\rightarrow {\cal O}_{\pq}(-2)^6\qfl{M}{\cal O}_{\pq}(-1)^6
\longrightarrow E\rightarrow 0$ (\ref{res})
gives by restriction to $Q$ a resolution
$$0\rightarrow {\cal O}_{Q}(-2)^6\qfl{M}{\cal O}_{Q}(-1)^6
\longrightarrow F\rightarrow 0\eqno(\ref{lag}.a)$$hence
$E_{|S}$ belongs to $H^{-1} (X)$. Conversely, any $F\in H^{-1} (X)$
admits a resolution (\ref{lag}.{\it a}), with $\pf M=0$  on $X$; the 
matrix $M$ defines an injective map ${\cal O}_{\pq}(-2)^6\rightarrow
{\cal O}_{\pq}(-1)^6$, whose cokernel is a bundle $E\in\mo_X$ such
that $E_{|S}\cong F$. Thus $r$ induces an isomorphism  $\mo_X\iso
H^{-1} (X)$. 
\ind Finally by (\ref{tyurin}) the smooth part of any fibre of $H$ is
Lagrangian.\cqfd
\rem{Remarks} a) The moduli space $\mo_S$ of stable rank 2 bundles on 
a K3 surface $S$ with
$c_1=0$, $c_2=4$ is a very interesting object (again, $c_2=4$ is the
smallest value for which the moduli space is not empty). O'Grady proved
that it is irreducible, and admits a smooth compactification
$\widehat{\mo}_S$ which is holomorphic symplectic [O]. It seems that 
$H$ does not extend to a Lagrangian fibration
$\widehat{\mo}_S\rightarrow \Pi$; however it is conceivable that it
extends to some other smooth birational model of $\widehat{\mo}_S$.
\ind b) A lagrangian fibration is the same thing as a completely 
integrable system: if 
$H=(H_1,\ldots,H_5)$ in some local coordinates system on $\Pi$, the
functions $H_i$ Poisson commute and the hamiltonian vector fields
$X_{H_i}$ span the tangent space to the fibre at each smooth point of
$H$. Since the fibres are open subsets of abelian varieties, it seems
likely that these vector fields  {\it linearize} on each fibre, that 
is, come from global vector fields on the abelian variety
(such a system is often called {\it algebraically} completely
integrable). To prove this it would be enough to exhibit some partial
compactification
$H':\mo'_S\rightarrow \Pi$ of $H$ such that $H'^{-1} (X)$, for $X$
general in $\Pi$, is isomorphic to the complement of a
subvariety of codimension $\ge 2$ in $J^2(X)$.

\vskip2cm
\centerline{ REFERENCES} \vglue15pt\baselineskip12.8pt
\def\num#1{\smallskip\item{\hbox to\parindent{\enskip [#1]\hfill}}}
\parindent=1.3cm 
 \num{B1}  A. {\pc BEAUVILLE}: {\sl Sous-vari\'et\'es sp\'eciales des
vari\'et\'es de Prym.}  Compositio math. {\bf 45}   (1981), 357--383.
\num{B2} A. {\pc BEAUVILLE}: {\sl 	Les singularit\'es du diviseur
$\Theta$ de la jacobienne interm\'ediaire de l'hypersurface cubique 
dans ${\bf P}^4$.} Algebraic threefolds, LN {\bf 947},
190--208; Springer-Verlag, 1982. 
\num{B3}  A. {\pc BEAUVILLE}: {\sl Determinantal hypersurfaces}. 
Michigan Math. Journal, to appear;
preprint math.AG/9910030.
\num{C-G} H. {\pc CLEMENS}, P. {\pc GRIFFITHS}: {\sl The intermediate 
Jacobian of the cubic threefold}. Ann. of Math.  {\bf 95}  (1972),
281--356.
\num{D} S. {\pc DRUEL}: {\sl Espace des modules des faisceaux 
semi-stables de rang $2$ et de classes de Chern $c_{1}=0$, $c_{2}=2$
et $c_{3}=0$ sur une hypersurface cubique lisse de}
${\bf P}^{4}$. Preprint math.AG/0002058.
\num{De} P. {\pc DELIGNE}: {\sl Le th\'eor\`eme
de Noether}.  SGA 7 II,  LN {\bf 340},  328--340;
Springer-Verlag,  1973.
\num{I} A. {\pc ILIEV}: {\sl Minimal sections of conic
bundles}. Boll. Unione Mat. Ital. {\bf 2}  (1999),
401--428. 
\num{I-M} A. {\pc ILIEV}, D. {\pc MARKUSHEVICH}:
{\sl  The Abel-Jacobi map for a cubic threefold and
periods of Fano threefolds of degree 14}. Preprint
math.AG/9910058.
\num{M-T} D. {\pc MARKUSHEVICH}, A. {\pc TIKHOMIROV}: {\sl  
The Abel-Jacobi map of a moduli component of vector bundles
on the cubic threefold}. J. of Algebraic Geometry, to
appear; preprint math.AG/9809140.
\num{Mo} B. {\pc MOISHEZON}: {\sl On $n$-dimensional compact
complex manifolds having $n$ algebraically independent meromorphic
functions} III. Amer. Math. Soc. Transl. (2) {\bf 63} (1967), 51--177.
\num{M} D. {\pc MUMFORD}: {\sl Lectures on curves on an
algebraic  surface}.  Annals of Math.
Studies {\bf 59}. Princeton University Press,  1966.
\num{Mu} S.  {\pc MUKAI}: {\sl Symplectic structure of the moduli
space of  sheaves on an abelian or $K3$ surface}. Invent. Math.
{\bf 77}  (1984), 101--116.
\num{O} K. O'{\pc GRADY}: {\sl Desingularized moduli spaces of sheaves
on a $K3$}. J. Reine Angew. Math. {\bf 512} (1999), 49--117.
\num{T1} A. N. {\pc TYURIN}: {\sl  The geometry of the Fano surface of
a nonsingular cubic $F\i\pq$ and Torelli theorems for Fano surfaces 
and cubics}. Math. USSR Izvestija {\bf 5} (1971), 517--546.
\num{T2} A. N. {\pc TYURIN}: {\sl Five lectures on three-dimensional 
varieties}. Russian Math. Surveys 27 (1972), no. 5, 1--53.
\num{T3} A. N. {\pc TYURIN}: {\sl The moduli spaces of vector bundles
on threefolds, surfaces and curves} I. Preprint, Erlangen, 1990.
\num{W} J. {\pc WAHL}: {\sl Gaussian maps on algebraic
curves}. J. Diff. Geometry {\bf 73} (1990), 77--98. 
\num{We} G. {\pc WELTERS}: {\sl Abel-Jacobi isogenies for
certain types of Fano threefolds}.  Math. Centre
Tracts, 141;  Math. Centrum, Amsterdam, 1981.

\vskip1cm
\def\pc#1{\eightrm#1\sixrm}
\hfill\vtop{\eightrm\hbox to 5cm{\hfill Arnaud {\pc BEAUVILLE}\hfill}
 \hbox to 5cm{\hfill DMA -- \'Ecole Normale
Sup\'erieure\hfill} \hbox to 5cm{\hfill (UMR 8853 du CNRS)\hfill}
\hbox to 5cm{\hfill  45 rue d'Ulm\hfill}
\hbox to 5cm{\hfill F-75230 {\pc PARIS} Cedex 05\hfill}}
\end